\documentclass{jT}
\usepackage{amssymb,latexsym}
\usepackage{amsmath}
\usepackage{amsthm}
\usepackage[english]{babel}
\DeclareMathOperator*{\vol}{vol}

\begin{document}
\thanks{The first author was
 partially supported by  NSF grant DMS-0401318, and PSC CUNY Research Award,
 No. 69288-00-38}

\newtheorem{introtheorem}{Theorem}
\renewcommand{\theintrotheorem}{\Alph{introtheorem}}
\newtheorem{theorem}{Theorem }[section]
\newtheorem{lemma}[theorem]{Lemma}
\newtheorem{corollary}[theorem]{Corollary}
\newtheorem{proposition}[theorem]{Proposition}
\theoremstyle{definition}
\newtheorem{definition}[theorem]{Definition}
\newtheorem{example}[theorem]{Example}
\newtheorem{remark}[theorem]{Remark}

\renewcommand{\labelenumi}{(\roman{enumi})} 
\def\theenumi{\roman{enumi}}

\numberwithin{equation}{section}

\def \g {{\gamma}}
\def \G {{\Gamma}}
\def \l {{\lambda}}
\def \a {{\alpha}}
\def \b {{\beta}}
\def \f {{\phi}}
\def \r {{\rho}}
\def \R {{\mathbb R}}
\def \H {{\mathbb H}}
\def \N {{\mathbb N}}
\def \C {{\mathbb C}}
\def \Z {{\mathbb Z}}
\def \F {{\Phi}}
\def \Q {{\mathbb Q}}
\def \e {{\epsilon }}
\def \ev {{\vec\epsilon}}
\def \ov {{\vec{0}}}
\def \GinfmodG {{\Gamma_{\!\!\infty}\!\!\setminus\!\Gamma}}
\def \GmodH {{\Gamma\backslash\H}}
\def \slr  {{\hbox{SL}_2( {\mathbb R})} }
\def \psl  {{\hbox{PSL}_2( {\mathbb R})} }
\def \L  {{\hbox{L}^2}}

\newcommand{\norm}[1]{\left\lVert #1 \right\rVert}
\newcommand{\abs}[1]{\left\lvert #1 \right\rvert}
\newcommand{\modsym}[2]{\left \langle #1,#2 \right\rangle}
\newcommand{\inprod}[2]{\left \langle #1,#2 \right\rangle}
\newcommand{\Nz}[1]{\left\lVert #1 \right\rVert_z}
\newcommand{\tr}[1]{\operatorname{tr}\left( #1 \right)}

\title[Hyperbolic lattice-point counting]{Hyperbolic lattice-point counting\\ and modular symbols}
\author{Yiannis N. Petridis}
\address{Department of Mathematics, University College London, Gower Street, London WC1E 6BT \\\newline The Graduate Center, Mathematics Ph.D. Program, 365 Fifth Avenue, Room 4208, New 
York, NY 10016-4309}
\email{petridis@member.ams.org}
\author{Morten S. Risager}
\address{Department of Mathematical Sciences, University of Aarhus, Ny Munkegade Building 530, 8000 {Aa}rhus, Denmark}
\email{risager@imf.au.dk}
\date{\today}


\maketitle
\begin{resume}
Soit un sous-group $\G$ de $\slr$ au quotient compacte et soit $\a$ une forme harmonique r\'eelle (non z\'ero). Nous \'etudions le comportement asymptotique de la fonction comptant des points du
r\'eseau hyperbolique $\G$ sous restrictions impos\'ees par des symboles modulaires $\modsym{\gamma}{\a}$. Nous montros que les valeurs normalis\'ees des symboles modulaires, ordonn\'ees selon ce comptage
poss\`edent une r\'epartition gausienne.
\end{resume}

\begin{abstr}
For a cocompact group $\G$ of $\slr$ we fix a real non-zero harmonic $1$-form
$\alpha$. We study the asymptotics of the hyperbolic lattice-counting problem for $\G$ under restrictions imposed by the modular symbols $\modsym{\gamma}{\a}$. We prove that the normalized values of the modular symbols, when ordered according to this counting, have a Gaussian distribution.
\end{abstr}

\section{Introduction}  
For a compact hyperbolic surface Huber \cite{huber} found the asymptotics in  the hyperbolic lattice counting problem. He also studied the distribution of the lengths of the closed geodesics. 
More precisely, let $X$ be a hyperbolic  surface. Let $\Gamma =\pi_1(X)$ be its fundamental group. 
Let $z$ and $w$ be two points in $\H$, the universal cover of $X$. Define the hyperbolic lattice counting function
\begin{equation*}N(z, w, x)=\{\gamma  \in \Gamma, r( \gamma z, w)\le x\},\end{equation*}
 where $r(z_1, z_2)$
denotes the hyperbolic distance between $z_1$ and $z_2$ in the $\H$. 
Huber \cite{huber} proved that
\begin{equation}\label{huberresult}N(z, w, x)\sim \frac{\pi}{\vol (X)}e^x.\end{equation}
One may also obtain error terms. We do not investigate these but
refer to \cite[Theorem 12.1]{iwaniec}, \cite{patterson,good, phillips-rudnick,chamizo}.
To every conjugacy class $\{\gamma\}$ of $\Gamma$ corresponds
a unique closed oriented geodesic of length $l(\gamma)$. Let $\pi (x)=\#\{\{\gamma\}\vert
 \textrm{$\g$ prime},  l(\gamma)\le x\}$.
The prime number theorem for closed geodesics states that \begin{equation}\pi (x)\sim e^x/x\end{equation} 
as $ x\rightarrow \infty$ and can be proved using the Selberg trace formula \cite{iwaniec}.
Generalizations of the hyperbolic lattice counting problem to infinite volume groups have been obtained by Lax and Phillips \cite{lax}. Recent applications of the hyperbolic lattice counting problem in this setting include sieving for almost primes among e.g. sums of squares, see  \cite{sarnak} and  \cite{kontorovich}.

On the other hand we have the  Poincar\'e pairing between homology and
cohomology: 
\begin{equation*}H_1(X,\R)\times H^1_{\rm dR}(X,\mathbb R)\rightarrow \mathbb R\end{equation*} and a projection
 $\phi:\Gamma\rightarrow H_1(X,\Z)$. Let $\modsym{\cdot}{\cdot}$ be the
composition of the two maps:
$$\modsym{\gamma}{\a}=\int_{\phi (\gamma)}\a.$$
We fix a real 
$1$-form $\a$. The numbers $\modsym{\gamma}{\a}$ will be called modular symbols, even in the compact setting. See \cite{goldfeld1, goldfeld2} for their arithmetic significance in the cofinite case of congruence subgroups.  In previous work the authors \cite{petr,  petr-ris, risager} have studied the distribution of the normalized values of the Poincar\'e pairing for compact and finite volume hyperbolic surfaces. In all articles we found as limiting distribution the normal Gaussian distribution. In all cases considered we were forced to group together elements in the same coset for appropriate coset decomposition of $\G$. In \cite{petr} we ordered the group elements of a cofinite subgroup $\G$ by realizing 
$\G=\pi_1(X)$ as a discrete subgroup of $\slr$, setting $\gamma=\left(\begin{array}{cc}a &b\\c&d\end{array}\right)$ and ordering the cosets $\gamma\in \GinfmodG$ according to $c^2+d^2$. To find the asymptotics of the powers of the modular symbols in \cite{petr} we used  Eisenstein series twisted by modular symbols \cite{goldfeld1, goldfeld2, gos}  which are sums over the cosets of $\GinfmodG$.   In \cite{risager} the matrix elements are ordered according to $(a^2+b^2)(c^2+d^2)$, and grouped together  in the same coset of a hyperbolic subgroup $\Gamma_{\gamma_1}$. Moreover, Risager defined in \cite{risager} and used hyperbolic Eisenstein series twisted by modular symbols to study the appropriate moments.

 In \cite{petr-ris} we considered the distribution of the values of the Poincar\'e pairing
where we order the conjugacy classes  of $\Gamma$ according to the lengths $l(\gamma)$. We proved that, if 
$$[[\gamma, \a]]=\sqrt{\frac{\vol(X)}{2\norm{\a}^2l(\gamma)}}\int_{\phi(\gamma )}\a.$$
 Then
 $$\frac{\#\left\{\gamma\in \pi_1(X)\vert [[\gamma, \a]]\in [a, b],
     l(\gamma)\le x\right\}}{\#\{\gamma\in\pi_1(X)\vert l(\gamma )\le x\}}\rightarrow\frac{1}{\sqrt{2\pi}}\int_a^be^{-t^2/2}\, dt $$
 as $x\rightarrow\infty$.
 To prove this theorem  we used the Selberg trace formula, via the Selberg zeta function and its derivatives
 in character varieties.
  
 In this paper we show how we can calculate the moments of the modular symbols in another ordering of the group elements $\gamma\in \Gamma$ \emph{but without using any coset decomposition}. This is the same ordering used by Huber in the hyperbolic lattice counting problem. Moreover, we have two parameters to play: the points $z$ and $w$ in the upper half plane $\H$.  The results are as follows.
 \begin{theorem}\label{maintheorem} Fix $\a$ a real harmonic $1$-form on $X$ of norm $\norm{\a}$  and two points $z, w\in \H$. 
 Let $$[\gamma, \a]=\sqrt{\frac{\vol(X)}{2\norm{\a}^2 r(\gamma z, w)}}\int_{\phi(\gamma )}\a.$$
 Then
 $$\frac{\#\left\{\gamma\in \Gamma \vert [\gamma, \a]\in [a, b], r(\gamma z, w)\le x\right\}}{\#\{\gamma\in \Gamma \vert r(\gamma z, w )\le x\}}\rightarrow\frac{1}{\sqrt{2\pi}}\int_a^be^{-t^2/2}\, dt $$
 as $x\rightarrow\infty$.
 \end{theorem}
The normalization depends on $z$, $w$ (but this is suppressed in the notation). However, the limiting distribution is universal. 
Along the way of proving  Theorem \ref{maintheorem} we introduce a new type of Dirichlet series twisted by modular symbols. Their study and analytic properties are of independent interest. Here are the results:
\begin{theorem}\label{anothertheorem}
Let 
\begin{equation*}G^{(n)}(z, w, s, 0)=\sum_{\gamma \in \Gamma}\frac{(-i)^n\modsym{\gamma}{\a}^n}{(\cosh(r(\gamma z, w))^s}.\end{equation*}
For each $n\in {\mathbf N}\cup \{0\}$ this series converges absolutely for $\Re (s)>1$. It admits analytic continuation in the whole complex plane with poles included in the set 
$$P=\{-2m+s_j \vert m \in {\mathbf N}\cup \{0\}, s_j(1-s_j)\in \rm{Spec} (\Delta)\}.$$
The order of the pole at $s=1$ is less than or equal to $[n/2]+1$, with order exactly $[n/2]+1$, for $n$ even. The poles at the other points of $P$ are of order less than equal to $n+1$. The functions $G^{(n)}(z, w, s, 0)$  are also bounded polynomially on vertical lines with $\Re (s)>1/2$.
\end{theorem}
Using Theorem \ref{anothertheorem} and the fact that we can find explicit
expressions for the leading term in the polar expansions when $n$ is even we can calculate all the asymptotic moments of the (normalized) modular symbols. Theorem \ref{maintheorem} then follows, see Section \ref{Sectionmoments}.

In principle we can find also the leading term in the polar expansion
when $n$ is odd, which gives asymptotic moments of the non-normalized
modular symbols.
As an example the first moment is given by the following: 
\begin{equation}\label{firstmoment}\sum_{\substack{\g\in \G\\ r(\gamma z, w)\leq x}}\modsym{\g}{\a}=
-\int_w^z \a \,\frac{\pi}{\vol (X)}e^x.\end{equation}
\section{Dirichlet series}
Let $\G$ be a discrete subgroup of $\psl$ with compact quotient $X=\GmodH$. Here $\H$ is the upper half-plane. Let $\chi$ be a unitary character of $\Gamma$. 
For our purposes we need in fact a
family of characters 
$$\begin{array}{lccc}\chi(\cdot,\e):&\Gamma&\rightarrow& S^1\\
                   &\g&\mapsto &\exp(-i\e\modsym{\g}{\alpha}).
\end{array}
$$
We have the following lemma: 
\begin{lemma}\label{eichlerbound}
\begin{equation*}
\modsym{\g}{\a} =\int_{\phi (\gamma )}\a=O(\norm{\a}_{\infty} r(\gamma z, w)),
\end{equation*}
where the implied constant depends on $z, w$.
\end{lemma}
\begin{proof}
By well-known properties of line integrals
$$\abs{\int_{\phi (\gamma )}\alpha }\le \norm{\alpha}_{\infty}l(\gamma ).$$
Since  $l(\gamma ) =\inf_{x\in X} r(x, \gamma x)$, we have
\begin{align*}\abs{\int_{\phi (\gamma )}\alpha }&\le
  \norm{\alpha}_{\infty}r(z, \gamma z)\\ &\le \norm{\alpha}_{\infty}(r(z, w)+r(w, \gamma z))\le \norm{\alpha}_{\infty}(D+r(\gamma z, w)), \end{align*}
 where $D$ is the diameter of $X$.
\end{proof}

Huber \cite{huber} introduced the series
$$G(z, w, s)=\sum_{\gamma \in \Gamma}\frac{1}{(\cosh r(\gamma z, w))^s},$$
which converges absolutely for $\Re (s)>1$.  Moreover, see \cite[Satz 1]{huber}, he proved that for fixed $w$ it is a continuous and automorphic function of $z$. It is twice differentiable in $z$.
We introduce the series
\begin{equation}\label{Gseries}G(z, w, s, \e)=\sum_{\gamma\in \Gamma}\frac{\chi (\gamma, \e)}{(\cosh r(\gamma z, w))^s}\end{equation}
and its $n$-th derivative with respect to $\e$, evaluated at $\e=0$, namely
\begin{equation}\label{Gnseries}G^{(n)}(z, w, s, 0)=\sum_{\gamma \in \Gamma}\frac{(-i)^n\modsym{\gamma}{\a}^n}{(\cosh(r(\gamma z, w))^s}.\end{equation}
The series $G(z, w, s, \e)$ converges absolutely for $\sigma=\Re (s)>1$ by comparison with $G(z, w, \sigma)$, as the character $\chi (\cdot, \e)$ is unitary. The series $G^{(n)}(z, w, s, 0)$ converges absolutely in the same region, using Lemma \ref{eichlerbound}.
 Since $G(z, w, s, \e)$ are automorphizations of $(\cosh r)^{-s}$,  they satisfy the shifted eigenvalue equation
\begin{equation}\label{shiftedequation}\Delta G(z, w, s, \e)+s(1-s)G(z, w, s, \e )=-s(s+1)G(z, w, s+2, \e).
\end{equation}
Here $$\Delta=y^2\left(\frac{\partial^2}{\partial
    x^2}+\frac{\partial^2}{\partial y^2}\right)$$ is the Laplace
operator for the upper half-plane, acting in the $z$ variable. The equation (\ref{shiftedequation}) appeared already in \cite{huber}.
The series (\ref{Gseries}), (\ref{Gnseries}) are analogous to the series 
\begin{equation*}\label{smaglos}
E(s, \e)=\sum_{\{\gamma_0\}}\frac{\chi(\g_0,\e)\ln (N(\gamma_0))}{N(\gamma_0)^s}
\end{equation*}
and its derivatives at $\e=0$
\begin{equation*}\label{marcopolo}
E^{(n)}(s,0)=\sum_{\{\gamma_0\}}\frac{{(-i)^n}\modsym{\g_0}{\a}^n\ln (N(\gamma_0))}{N(\gamma_0)^s}.
\end{equation*}
studied in our previous work \cite{petr-ris}. Here the  sums are over only primitive conjugacy classes $\{\gamma_0\}$ of $\Gamma$.

We need to show that (\ref{Gnseries}) admits meromorphic continuation
to the whole complex plane, find the pole order at $s=1$  (Section
\ref{sectionpole}) and find the leading term in the Laurent expansion of $G^{(n)}(z, w, s, 0)$ (Theorem \ref{poles}).  We also need to control these derivatives as functions of $s$ on vertical lines for $\Re (s)>1/2$ to be able to use contour integration techniques 
from analytic number theory (Section \ref{sectiongrowth}).
Theorem \ref{anothertheorem} summarizes (less precisely)  these results.

We consider the space
$$\L(\GmodH,\overline \chi(\cdot,\epsilon))$$ of
$(\G,\overline \chi(\cdot,\epsilon))$-automorphic functions, i.e. functions
$f:\H\to\C$ satisfying $$f(\g z)=\overline{\chi}(\g,\e)f(z),$$ and
$$\int_{\GmodH}\abs{f(z)}^2d\mu(z)<\infty.$$  Here $d\mu(z)=y^{-2}dxdy$ is
the invariant Riemannian measure on $\H$ derived from the Poincar\'e
metric $ds^ 2=y^{-2}(dx^2+dy^2)$. We shall denote by $\norm{\cdot}$
the usual norm in the Hilbert space $\L(\GmodH,\overline
\chi(\cdot,\epsilon))$. It is easy to see that $G(z, w, s, \e)\in\L(\GmodH,\overline \chi(\cdot,\epsilon))$. The automorphic Laplacian $\tilde L(\e)$ is
the closure of the operator $\Delta$ acting on smooth functions in
$\L(\GmodH,\overline \chi(\cdot,\epsilon))$. The  spectrum of $\tilde{L}(\e)$ is discrete and $-\tilde L(\e)$ is nonnegative with eigenvalues
$$0\leq \l_0(\e)\le \l_1(\e)\leq\ldots .$$ As usual we set $\l_j(\e)=s_j(\e)(1-s_j(\e))$.
  The first eigenvalue is zero if and only if $\e=0$ and in this case it is a simple eigenvalue.

The resolvent $\tilde R(s,\e)=(\tilde L(\e)+s(1-s))^{-1}$, defined off the spectrum of $\tilde L(\e)$, is a Hilbert-Schmidt operator on $\L(\GmodH,\overline \chi(\cdot,\epsilon))$. It is a holomorphic (in $s$) family of operators, and the operator norm of the resolvent is bounded as follows:
\begin{equation}\label{calculustestsareboring}
\norm{\tilde R(s,\e)}_\infty\leq \frac{1}{\hbox{dist}({s(s-1),\hbox{spec}(\tilde L(\e))})}\leq \frac{1}{\abs{t}(2\sigma-1)},
\end{equation} where $s=\sigma+it$, $\sigma>1/2$.

We fix $z_0\in\H$ and introduce unitary operators (see (\cite{phillipssarnak2})
\begin{eqnarray}
U(\e):&\L(\GmodH)&\to \L (\GmodH,\overline \chi(\cdot,\epsilon))
\\
&f&\mapsto\exp\left(i\e\int_{z_0}^z\alpha\right)f(z). \nonumber 
\end{eqnarray}
We then define
\begin{eqnarray}
L(\e)&=&U^{-1}(\e)\tilde L(\e)U(\e),\\
R(s,\e)&=&U^{-1}(\e)\tilde R(s,\e)U(\e). 
\end{eqnarray}
This ensures that $L(\e)$ and $R(s,\e)$ act on the fixed space $\L(\GmodH)$. It is then easy to verify that 
\begin{align}\label{muffinstogo}L(\e)h=\Delta h +2i\e\modsym{dh}{\alpha}&-i\e\delta(\alpha)h-\e^ 2\modsym{\a}{\a}h,\\
\label{coffeetogo}(L(\e)+s(1-s))R(s,\e)=&R(s,\e)(L(\e)+s(1-s))=I.\nonumber
\end{align}
Here \begin{eqnarray*}\modsym{f_1dz+f_2d\overline z}{g_1dz+g_2d\overline z}&=&2y^ 2(f_1\overline{g_1}+f_2\overline{g_2})\\
\delta(pdx+qdy)&=&-y^2(p_x+q_y).
\end{eqnarray*}
We notice that $\delta(\alpha)=0$ since $\alpha$ is harmonic. We notice also that 
\begin{align}
L^{(1)}(\e)h&=2i\modsym{dh}{\a}-2\e\modsym{\a}{\a}h,\\ \allowdisplaybreaks
\label{second} L^{(2)}(\e)h&=-2\modsym{\a}{\a}h,\\\allowdisplaybreaks
L^{(i)}(\e)h&=0,\quad\textrm{ when }i\geq 3.
\end{align}
(We use superscript $(n)$ to denote the
$n$-th derivative in $\e$.) Using Lemma \ref{eichlerbound} we see that
$$\int_{z_0}^{\gamma z}\alpha=O(r(\gamma z, w)).$$
In order to investigate the behavior of $G^{(n)}(z,w,0)$ we define the
auxiliary series 
$$D(z, w, s, \e)=U(-\e) G(z, w, s, \e)=\sum_{\gamma\in \Gamma}\frac{\exp(-i\e\int_{z_0}^{\gamma z}\alpha)
}{(\cosh r(\gamma z, w ))^s}\in \L(\GmodH).$$
This series converges absolutely for $\Re (s)>1$.
Using (\ref{shiftedequation}) we see that
\begin{equation}\label{equationD}(L(\e)+s(1-s))D(z, w, s, \e)=-s(s+1)D(z, w, s+2, \e).\end{equation}
Since $D(z,w,s\e)$ is square integrable on $\GmodH$ it follows that\begin{equation}\label{LP}D(z, w, s, \e)=-s(s+1)R(s, \e)D(z, w, s+2, \e),\end{equation}
and this implies the meromorphic continuation of $D(z, w, s, \e)$ in the whole complex plane with (at most) poles of first order at the spectral points $s_j$ and $-2n+s_j$, $n=0, 1, \ldots$. This is done in the standard way of extending the domain of analyticity on subsequent vertical strips of width $2$ to the left, using the well-known meromorphic continuation of the resolvent (with poles of order $1$ at the points $s_j$).
The series 
$$D^{(n)}(z, w, s, 0)=\sum_{\gamma\in \G}\frac{(-i)^n(\int_{z_0}^{\g z}\alpha )^n}{(\cosh r(\g z, w))^s}$$
converges absolutely for $\Re (s)>1$. 
By differentiating (\ref{equationD}) $n$-times and plugging $\e=0$ we get
\begin{align}\nonumber(\Delta+s(1-s))D^{(n)}(z, w, s,
  0)& \\\nonumber +\binom{n}{1}L^{(1)}D^{(n-1)}&(z, w, s, 0) +\binom{n}{2}L^{(2)}D^{(n-2)}(z, w, s, 0)\\\label{knut}=&-s(s+1)D^{(n)}(z, w, s+2, 0).\end{align}
This implies for $n\ge 2$
\begin{align}\label{FF}\nonumber{D^{(n)}(z, w, s, 0)}=-s(s+1)&R(s,
  0)D^{(n)}(z, w, s+2, 0)\\ -\binom{n}{1}&R(s, 0)L^{(1)}D^{(n-1)}(z,
  w, s, 0)\\ \nonumber & -\binom{n}{2}R(s, 0)L^{(2)}D^{(n-2)}(z, w, s, 0), \end{align}
while 
\begin{eqnarray}\label{tralala}D^{(1)}(z, w, s, 0)&=&-s(s+1)R(s, 0)D^{(1)}(z, w, s+2, 0)\\&&-R(s, 0)L^{(1)}D(z, w, s, 0)\nonumber. \end{eqnarray}
These equations prove the meromorphic continuation of $D^{(n)}(z, w, s, 0)$ in the whole complex plane with poles at $-2n+s_j$, $n =0, 1, \ldots$ using induction on $n$ and continuation in vertical strips from $\Re (s)>1$ to $-1<\Re (s)\le 1$, and subsequently to $-3<\Re (s)\le -1$, etc. It also follows that $D^{(n)}(z, w, s, 0)$ is analytic for $\Re (s)>1$, and that at points where is it analytic it is square integrable.
We also see that 
\begin{equation}\label{transport}G^{(n)}(z, w, s, 0)=\sum_{j=0}^n\binom{n}{j}\left(i\int_{z_0}^z\a\right)^jD^{(n-j)}(z, w, s, 0).\end{equation}

\section{The pole at $s=1$.}\label{sectionpole}
In this section we identify the pole order and the leading term of the
function
$G^{(n)}(z, w, s,0)$ at $s=1$. In $G^{(n)}(z, w, s,\e)$, $R(s,\e)$ and
$L^{(n)}(\e)$ we shall often omit $0$ from the notation when we
set $\e=0$.  We note that $G(z, w, s)$ has a first order pole with residue $2\pi/\vol({\GmodH})$, as shown by Huber \cite[Satz 2, p.~5]{huber}. We derive it below only for completeness. 

We recall that close to $s=1$ \begin{equation}\label{resolventexpansion}R(s)=\sum_{i=-1}^\infty R_i (s-1)^i,\qquad R_{-1}=-P_0\end{equation}  and that $R(s)-R_{-1}(s-1)^{-1}$ is holomorphic in $\Re(s)>h$. Here  $h=\Re(s_1)$, $s_1(1-s_1)=\l_1$ is the first non-zero eigenvalue, and \begin{equation*}P_0f=\langle f,\vol(\GmodH)^{-1/2}\rangle\vol{(\GmodH)}^{-1/2}\end{equation*} is the projection of $f$ to the zero eigenspace.

For $n=0$ the residue of $D(z, w, s)$ at $s=1$ is from (\ref{LP})
\begin{align*}\frac{-1(2)(-1)}{\vol(\GmodH)}\int_{\GmodH}D(z, w,
  3)d\mu (z)&=\frac{2}{\vol(\GmodH)}\int_0^{\pi}\int_0^{\infty}
  \frac{2\sinh r}{(\cosh r)^3}\, dr\,d \phi\\ &=\frac{2\pi}{\vol(\GmodH)}.\end{align*}
Here we have unfolded the integral $\int_{\GmodH}$ to the whole upper half-plane and integrated in polar coordinates, where the invariant hyperbolic measure is $2\sinh r\,dr\, d\phi$, $r\ge0$, $\phi\in[0, \pi]$.
 
Close to $s=1$ we also have the expansion of $D(z, s, w)$ as
$$\frac{2\pi/\vol(\GmodH)}{s-1}+D_0(z, w)+D_1(z, w)(s-1)+\cdots .$$

The crucial observation is that 
\begin{equation}\label{everythingdies}
L^{(1)} P_0=0, \qquad P_0 L^{(1)} \subseteq 0.\end{equation}
The first equality follows from the fact that $L^{(1)}$ is a differentiation operator while $P_0$ projects to the constants. The second equality follows from the first by using that both operators are self-adjoint.

For $n=1$ we look at (\ref{tralala}). We notice that $R(s)L^{(1)}D(z, w, s)$ is regular at $s=1$: at $s=1$ the singular terms are
$$\frac{-P_0 L^{(1)} (2\pi/\vol(\GmodH))}{(s-1)^2}+\frac{-P_0 L^{(1)} D_0(z, w)+R_0 L^{(1)}( 2\pi/\vol(\GmodH))}{s-1}=0.$$

We can now prove:
\begin{theorem}\label{poles} For $n\geq 0$, $D^{(n)}(z, w, s)$ and $G^{(n)}(z, w, s)$ have a pole at $s=1$ of order at most $[n/2]+1$. If $n=2m$ the pole is of order $m+1$ and the $(m+1)$-term in the expansion around $s=1$ is $$(-1)^m\frac{(2\pi)(2m)! \norm{\a}^{2m}}{\vol(\GmodH)^{m+1}}\frac{1}{(s-1)^{m+1}}.$$
\end{theorem}
\begin{proof}

The claim has been shown above for $n=0$ and $n=1$. 
We assume that the order of the pole of $D^{(n-1)}(z, w, s)$ and $D^{(n-2)}(z, w, s)$ at $s=1$ are less than or equal to $[(n-1)/2]+1$ and $[(n-2)/2]+1$ respectively.  By (\ref{second}), (\ref{FF}), and (\ref{everythingdies}) the order of the pole of
  of $D^{(n)}(z, w, s)$ at $s=1$ is less than equal to $\max( [(n-1)/2]+1, [(n-2)/2]+1+1)=[n/2]+1$.
  
  Assume that for $n=2(m-1)$ we have proved the claim for the order of the pole and leading singularity. By (\ref{second}), (\ref{FF}) and (\ref{everythingdies}) the coefficient of $(s-1)^{-m-1}$ in $D^{(2m)}(z, w, s)$ is
  $$-\binom{2m}{2}\int_{\GmodH}\frac{-1}{\vol(\GmodH)} (-2\langle \a, \a \rangle)\frac{(-1)^{m-1}(2\pi) (2m-2)!\norm{\a}^{2m-2}}{\vol (\GmodH)^{m}}\, d\mu (z),$$
  and this is $$\frac{(2m)!}{2(2m-2)!}2\norm{\a}^2  \frac{(-1)^m(2\pi) (2m-2)!\norm{\a}^{2m-2}}{\vol (\GmodH)^{m+1}}=\frac{(-1)^m(2\pi)(2m)!\norm{\a}^{2m}}{\vol(\GmodH)^{m+1}}.$$
  The claim for $G^{(n)}(z, w, s)$ follows from (\ref{transport}) and the result for $D^{(n)}(z, w, s)$.
\end{proof}

\section{Growth on vertical lines}\label{sectiongrowth}
\begin{lemma} For $\sigma =\Re (s)\ge \sigma_0>1$ , we have $G^{(n)}(z, w, \sigma+ it)=O(1)$.\end{lemma}
This follows from the absolute convergence of the series in this region.

\begin{lemma}\label{ohweh} For $\sigma=\Re (s)>1/2$ we have $G(z, w, s)=O(\abs{t}^{6(1-\sigma)+\e}).$
\end{lemma}
\begin{proof}By (\ref{shiftedequation}) and (\ref{calculustestsareboring}) we get
$$\norm{G(z, w, s )}_2\le\frac{\abs{s(s+1)}}{\abs{t}(2\sigma -1)}\norm{G(z, w, s+2)}_2\ll \abs{t}$$ and
$$\norm{\Delta G(z, w, s )}_2\ll \abs{t}^3.$$
By the Sobolev embedding theorem, which provides the bound
\begin{equation}\label{sobolev}\norm{f}_{\infty}\ll \norm{f}_2+\norm{\Delta f}_2\end{equation}
we get
$$\norm{G(z, w, s )}_{\infty}=O(\abs{t}^3).$$
We apply the Phragm\'en-Lindel\"of principle in the strip $1/2+\e\le \sigma \le 1+\e$ to get the result.
\end{proof}
The right-hand side in (\ref{sobolev}) is the $H^2$-norm of $f$.
Similarly one can define the $H^1$-norm using any first order differential operator and we have $\norm{f}_{H^1}\ll \norm{f}_{H^2}.$
\begin{lemma}\label{growth} For $\Re (s)>1/2$ we have $$D^{(n)}(z, w, s)=O(\abs{t}^{6(n+1)(1-\sigma)+\e}), \quad G^{(n)}(z, w, s)=O(\abs{t}^{6(n+1)(1-\sigma)+\e}).$$
\end{lemma}
\begin{proof} Since the series of $D^{(n)}(z, w, s)$ and $G^{(n)}(z, w, s)$ converge absolutely for $\Re (s)>1$, we get that they are $O(1)$ in this region. To prove the lemma one works as in Lemma \ref{ohweh}, by first proving 
\begin{equation}\label{exponents}
 D^{(n)}(z, w, s)=O(\abs{t}^{3(n+1)}), \quad  G^{(n)}(z, w, s)= O(\abs{t}^{3(n+1)})
 \end{equation}
 for $\Re (s)>1/2$ and then use the Phragm\'en-Lindel\"of principle.
 We prove inductively the estimates
 \begin{eqnarray}\label{1} \norm{D^{(n)}(z, w, s)}_{\infty}&=&O(\abs{t}^{3(n+1)}), \\
 \label{2}\norm{L^{(1)}D^{(n-1)}(z, w, \sigma+it)}_2&=&O(\abs{t}^{3n+2}), \\ 
 \label{3}\norm{L^{(2)}D^{(n-2)}(z, w, s)}_2&=&O(\abs{t}^{3(n-1)}).\end{eqnarray}
Assume that the result is true for $m\le n-1$. We prove it for $n$.
Using (\ref{1}) for $n-2$ we have \begin{displaymath}\norm{D^{(n-2)}(z, w, s)}_{\infty}=O(\abs{t}^{3(n-1)})\Longrightarrow
\norm{L^{(2)}D^{(n-2)}(z, w, s)}_2=O(\abs{t}^{3(n-1)})
,\end{displaymath}
as $L^{(2)}$ is a multiplication operator on a compact set. This proves (\ref{3}) for $n$.
We have
\begin{align*}
\norm{L^{(1)}D^{(n-1)}(z, w, s)}_2&\le c_1\norm{D^{(n-1)}(z, w, s)}_{H^1}\\
&\hspace{-1.5cm}\le c_2\norm{D^{(n-1)}(z, w, s)}_{H^2}\\
&\hspace{-1.5cm}\le c_2
\left(\norm{D^{(n-1)}(z, w, s)}_2+\norm{\Delta D^{(n-1)}(z, w, s )}_2\right)\\
&\hspace{-1.5cm}\ll\abs{t}^{3n}+\abs{t}^2\norm{D^{(n-1)}}_2+\norm{L^{(1)}D^{(n-2)}}_2+\norm{L^{(2)}D^{(n-3)}}_2\\
&\hspace{-1.5cm}\ll\abs{t}^{3n}+\abs{t}^{3n+2}+\abs{t}^{3n-1}+\abs{t}^{3(n-2)}=O(\abs{t}^{3n+2}).
\end{align*}
Here we have used the inductive hypothesis and (\ref{knut}).
This proves (\ref{2}) for $n$.
We have by (\ref{FF}) and (\ref{calculustestsareboring}) and the inductive hypothesis
\begin{equation}\label{flocke}
\norm{D^{(n)}(z, w, s)}_2\ll \frac{1}{\abs{t}}(\abs{t}^{3n+2}+\abs{t}^{3(n-1)})=O(\abs{t}^{ 
3n+1}).\end{equation}
Using (\ref{knut}) and 
(\ref{flocke}) and the previous results we get
\begin{align}\label{baer}
\nonumber\norm{\Delta D^{(n)}(z, w, s )}_2 \ll \abs{t}^2 &\norm{D^{(n)}(z, w, s
  )}_2+\norm{L^{(1)}D^{(n-1)}}_2+\norm{L^{(2)}D^{(n-2)}}_2\\  &\ll\abs{t}^{3n+3}
.\end{align}
Equations (\ref{flocke}) and  (\ref{baer})  prove (\ref{1}) for $n$.
 \end{proof}
\section{Calculating the moments}\label{Sectionmoments}
We are now ready to prove Theorem \ref{maintheorem}. The proof uses the method of asymptotic moments precisely as in \cite{petr,risager}. From Theorem \ref{poles}, Lemma \ref{growth} and Lemma \ref{eichlerbound} we may conclude, using a more or less standard contour integration argument (see  \cite{petr,risager} for details), that, as $T\to\infty$,
\begin{align}
\label{summatory} \sum_{\substack{\g\in \G\\ \cosh r(\gamma z, w)\leq
    T}}\!\!\!\!\!\!\!\!\!\!\!&\modsym{\g}{\a}^n\\ \nonumber &=\left\{\begin{array}{ll}\displaystyle{\frac{(2\pi)(2m)!\norm{\a}^{2m}}{m!\vol{(\GmodH)}^{m+1}}}T(\log T)^m+O(T(\log T)^{m-1}),&n=2m,\\
O(T(\log T)^m),&n=2m+1.
\end{array}\right.
\end{align} 
Setting $T=\cosh(x)$ and using $\cosh(x)=e^x/2+O(1)$ as $x\to \infty$
we can formulate this as follows:
\begin{equation}\label{summatory'}
\sum_{\substack{\g\in \G\\ r(\gamma z, w)\leq x}}\!\!\!\!\!\!\modsym{\g}{\a}^n=\left\{\begin{array}{ll}\displaystyle{\frac{\pi(2m)!\norm{\a}^{2m}}{m!\vol{(\GmodH)}^{m+1}}}e^xx^m+O(e^xx^{m-1}),&n=2m,\\
O(e^xx^m),&n=2m+1.
\end{array}\right.
\end{equation} 
as $x\to \infty$. Let now 
\begin{equation}
[\g,\a]=\sqrt{\frac{\vol({\GmodH})}{2 r(\g z, w)\norm{\a}^2}}\modsym{\g}{\a}.
\end{equation}
We then define the random variable $Y_x$ with probability measure 
\begin{equation}
P(Y_x\in [a,b])=\frac{
\#\{ \g\in \G  \vert r(\g z, w)\leq x,\quad[\g,\a] \in [a,b]\} }
{\#\{ \g \in \G \vert r(\g z, w)\leq x\}}.
\end{equation}
We want to calculate the asymptotic moments of these, i.e. find the limit of 
\begin{equation}M_n(Y_x)=\frac{1}{\#\{  \g \in\G \vert r(\g z, w)\leq x\}}\sum_{\substack{ \g\in \G\\ r(\g z, w)\leq x}}[\g ,\a]^n\end{equation} 
as $x\to\infty$. We note that by (\ref{huberresult}) or, alternatively, (\ref{summatory}), the denominator is asymptotically $\pi e^x/\vol (\GmodH)$. By
partial summation we have 
\begin{align*}\sum_{r(\g z, w)\leq x}&[\g,\a]^n=\frac{\vol(\GmodH)^{n/2}}{\norm{\a}^n2^{n/2}}\sum_{r(\g z, w)\leq x}\modsym{\g }{\a}^n\frac{1}{r(\g z, w)^{n/2}}\\
=&\frac{\vol(\GmodH)^{n/2}}{\norm{\a}^n2^{n/2}x^{n/2}}\sum_{r(\g z, w)\leq x}\modsym{\g }{\a}^n+O(e^x x^{-1}),
\end{align*}
This may be evaluated by (\ref{summatory'}). We find
$$M_n{(Y_x)}\to \begin{cases}\displaystyle{\frac{(2m)!}{m!2^m}},&\textrm{if
      $n=2m$,}\\
0,&\textrm{otherwise.}
\end{cases}$$
We notice that the right-hand side coincides with the moments of the
Gaussian distribution. Hence by a classical result due to Fr\'echet
and Shohat, see \cite[11.4.C]{loeve},  we may conclude that 
$$P(Y_x \in
[a,b])\to\frac{1}{\sqrt{2\pi}}\int_a^b\exp\left(-\frac{t^2}{2}\right)dt\textrm{ as $x\to\infty$}.$$
This concludes the proof of Theorem \ref{maintheorem}.

\begin{remark}
Equation (\ref{firstmoment}) follows from (\ref{tralala}) and (\ref{transport}) as follows. We take $z_0=w$. The pole of $$\left(i\int_w^z \a\right) D^{(0)}(z, w, s, 0)$$ contributes a  residue 
$$i\int_w^z \a \frac{2\pi}{\vol (\GmodH )}.$$
One the other hand the term $-s(s+1)R(s, 0)D^{(1)}(z, w, s, 0)$
contributes a residue $0$, since 
$$\int_{\GmodH}D^{(1)}(z, w, 3, 0)\, d\mu (z)=\int_{\H}\frac{-i\int_w^z \a}{\cosh r(z, w)^3}\, d\mu (z).$$
Now an integration in polar coordinates centered at $w$  produces the average of $F(z)=\int_w^z \a$ on a circle around $w$ and this is the average of a harmonic function on a circle. The mean-value theorem gives that it is equal to $2\pi F(w)=0$.
\end{remark}

\end{document}